\documentstyle[11pt,graphicx]{article}

\newcommand{\la}{\lambda}
\newcommand{\om}{\omega}

\newcommand{\Om}{\Omega}

\begin{document}
\title{{\bf Isoperimetric Inequalities and Sharp Estimate for Positive Solution of Sublinear Elliptic Equations}
   % Enter your title between curly braces
\thanks{Project supported by Natural Science Foundation of China (No. 10971061)}\thanks{Corresponding authors: Qiuyi Dai and Huaxiang Hu} }
\author{{\small  Qiuyi\hspace*{0.1cm}Dai\thanks{daiqiuyi@yahoo.com.cn} \hspace*{1.0cm}
Renchu\hspace*{0.1cm}He \hspace*{1.0cm}
Huaxiang\hspace*{0.1cm}Hu\thanks{hunanhhx@163.com}\hspace*{0.2cm}}\\
{\small Department of Mathematics, Hunan Normal University,
Changsha, Hunan, 410081, China}\\
}
        % Enter your name between curly braces
\date{}          % Enter your date or \today between curly braces
\maketitle

\small\begin{quote}\bf Abstract:\quad \rm In this paper, we prove
some isoperimetric inequalities and give a sharp bound for the
positive solution of sublinear elliptic equations.

{\bf Key words:}{\rm  \hspace*{0.05cm} Isoperimetric inequality,
Schwarz symmetrization, positive solution, sublinear equation.}

%\noindent\bf AMS classification:\quad \rm 4QJ10, 35P15, 49J20.
\end{quote}\normalsize

\section{Introduction and Main Results}

\setcounter{section}{1}
\renewcommand{\theequation}{\thesection.\arabic{equation}}

\vskip 0.05in

\noindent Let $\Omega\subset R^n$ be a bounded domain whose boundary
$\partial \Om$ is assumed to be of Lipschitz type. Assume that
$0<q<1$. We consider the following problem.
\begin{eqnarray}
\left
 \{\begin{array}{ll}
  -\Delta u=u^q, &x\in \Om,\\
   u>0,& x\in \Om\\
  u=0,& x\in\partial \Om.
\end{array}
\right.
\label{4.1.6}
\end{eqnarray}

\vskip 0.1cm

The purpose of this paper is to prove some isoperimetric
inequalities and give sharp bound for the solution of problem
(\ref{4.1.6}) by making use of rearrangement method.

\vskip 0.1cm

There are a lot of materials on isoperimetric inequalities for
eigenvalues and eigenfunctions of elliptic operators. For the
isoperimetric inequalities on eigenvalues of elliptic operators we
refer to \cite{ref CB, ref MSA, ref QD, ref BV, ref Bha, ref Bos,
ref Dan, ref Fab, ref Kra, ref Kra1, ref Sze, ref Wein1} and on
eigenfunctions we refer to \cite{ref PR1, ref PR2, ref GC, ref FC,
ref MVB, ref KJ1, ref KJ2}. The first result on isoperimetric
inequality for eigenfunctions of Laplace operator was obtained by
Payne and Rayner in \cite{ref PR1} . In 1972, Payne and Rayner
considered in \cite{ref PR1} the following eigenvalue problem
defined on bounded domains in $R^2$
\begin{eqnarray}
\left \{\begin{array}{ll}
-\Delta \varphi=\lambda\varphi & \mbox{in}\ \Om,\\
\varphi=0 & \mbox{on}\ \partial \Om,
\end{array}
\right.
 \label{4.1.2}
\end{eqnarray}
and prove that for the first eigenvalue $\lambda_1(\Omega)$ and the
first eigenfunction $\varphi_1(x)$ of problem (\ref{4.1.2}), the
following inequality holds
\begin{eqnarray}
\left(\int_{\Om}|\varphi_1|dA\right)^2\geq\frac{4\pi}{\lambda_1(\Omega)}\int_{\Om}\varphi_1^2dA,
\label{4.1.3}
\end{eqnarray}
with equality if and only if $\Om$ is a disk.

\vskip 0.1cm

Unfortunately, the argument used by Payne and Rayner works only for
the case $n=2$. Kohnler-Jobin \cite{ref KJ1, ref KJ2} and G.Chiti
\cite{ref GC} generalized the Payne and Rayner's inequality
(\ref{4.1.3}) to arbitrary dimension $n$ by employing the Schwarz
symmetrization method. It is by now well known that the Schwarz
symmetrization method is very useful for the estimate of sharp bound
of solutions to elliptic and parabolic equations, and has been
extensively studied since the pioneer works of Weinberger \cite{ref
Wein}, Talenti \cite{ref Tal} and Bandle \cite{ref Ban1}. See for
example \cite{ref GT, ref BM, ref ATL1, ref ATL2} for more details.
The basic idea in the use of the symmetrization method is to compare
the orignal problem with an auxiliary problem defined on a suitable
ball. Let $\Om^*$ be the Schwarz symmetrization of $\Omega$, that is
$\Omega^*$ is a ball in $R^n$ with center at $0$ and such that
$\Omega^*$ and $\Omega$ have same volume. The auxiliary problem used
by Kohnler-Jobin to generalize inequality (\ref{4.1.3}) can be read
as

\begin{equation}
\left\{\begin{array}{ll}\label{eq1}
-\Delta \varphi=\alpha\varphi+1 & \mbox{in}\ \Om^*,\\
\varphi=0 & \mbox{on}\ \partial \Om^*,
\end{array}
\right.
\end{equation}
with $-\infty<\alpha<\lambda_1(\Omega)$.

Whereas, G.Chiti used an auxiliary problem defined on  a ball
smaller than $\Om^*$ which can be read as
 \begin{eqnarray}
\left
 \{\begin{array}{ll}
  -\Delta z=\lambda z, &x\in B_r(0),\\
  z=0,& x\in\partial B_r(0),
\end{array}
\right.
 \label{5}
\end{eqnarray}
where $r=\sqrt{\frac{\lambda_1(\Om^*)}{\lambda_1(\Om)}}R^*$ and
$R^*$ is the radius of $\Om^*$.

It follows from the famous Faber-Krahn inequality that $r\leq R^*$,
and hence $B_r(0)$ is smaller than $\Om^*$. Furthermore, an easy
computation implies that the first eigenvalue of problem (\ref{5})
is $\lambda_1(\Om)$.

Compared with the auxiliary problem used by Kohnler-Jobin, the
problem used by  G.Chiti is more natural and extendable for other
situations.

Let $\varphi_1(x)$ be the first eigenfunction of problem
(\ref{4.1.2}), and $z_1(x)$ be the first eigenfunction of problem
(\ref{5}). If we normalize $\varphi_1(x)$ and $z_1(x)$ so that
$\int_{\Om}\varphi^p_1(x)dx=\int_{B_r(0)}z^p_1(x)dx$ for $p>1$, then
a celebrate result established by G.Chiti in \cite{ref GC} can be
stated as

\vskip 0.1in

{\bf Conclusion A}.\ There exists an unique point $s_0\in (0,\
|B_r(0)|)$ such that
 \begin{eqnarray*}
\left\{\begin{array}{ll}
z^*_1(s)>\varphi^*_1(s), &\mbox{for}\ s\in (0,\ s_0),\\
z^*_1(s)<\varphi^*_1(s), &\mbox{for}\ s\in (s_0,\ |B_r(0)|).
\end{array}
\right.
\end{eqnarray*}
where $ z^*_1(s)$ and $\varphi^*_1(s)$ are the decreasing
rearrangement of $z_1(x)$ and $\varphi_1(x)$ respectively, and
$|B_r(0)|$ denotes the volume of $B_r(0)$.

\vskip 0.1in

By making use of conclusion A, Chiti proved a reverse Holder
inequality for the first eigenfunction of problem (\ref{4.1.2})
which, in turn, is an isoperimetric inequality and more stronger
than inequality (\ref{4.1.3}). It is worth pointing out that a most
important application of conclusion A can be found in the proof of
the famous P.P.W conjecture (see \cite{ref MSA}).

Contrast to the eigenvalue problem, there are few results on the
isoperimetric inequalities for solutions of semilinear elliptic
problem. This is the motivation of our study of the isoperimetric
inequalities for the solution of problem (\ref{4.1.6}). Our method
is adapted from G. Chiti's paper \cite{ref GC} by carefully choosing
the comparison problem.

To state our results, we introduce the following auxiliary problem
 \begin{eqnarray}
\left
 \{\begin{array}{ll}
 -\Delta h=h^q, &x\in \Om^*,\\
   h>0,& x\in \Om^*\\
  h=0,& x\in\partial \Om^*.
\end{array}
\right.
 \label{4.1.7}
\end{eqnarray}
where $\Om^{*}$ is the Schwarz symmetrization of $\Omega$.

Let $\sigma_1=\frac{2(1+q)k+(1-q^2)n}{n+2-(n-2)q}$ and
$\sigma_2=\frac{2(1+q)}{n+2-(n-2)q}$ be fixed. Then our main result
can be stated as

\vskip 0.1in

{\bf Theorem 1.1.}  Let $u(x)$ be the unique solution of problem
(\ref{4.1.6}) and $h(x)$ be the unique solution of problem
(\ref{4.1.7}). Then for any $k\geq q+1$, we have
\begin{eqnarray}
\int_{\Om}u^k(x)dx\leq
C(q,k,\Om^{*})\|u\|^{\sigma_1}_{L^{q+1}(\Om)}.
\label{4.1.4}
\end{eqnarray}
Consequentely
\begin{eqnarray}
\max\limits_{x\in\Om}u(x)\leq C(q,\Om^*)\|u\|^{\sigma_2}_{L^{q+1}(\Om)},\label{4.1.5}
\end{eqnarray}
where
$C(q,k,\Om^{*})=\int_{\Om^{*}}h^k(x)dx/\|h\|^{\sigma_1}_{L^{q+1}(\Om^{*})}$
and
$C(q,\Om^*)=\max\limits_{x\in\Om^*}h(x)/\|h\|^{\sigma_2}_{L^{q+1}(\Om^{*})}$.
Moreover, the equality holds in each of inequalities (\ref{4.1.4})
and (\ref{4.1.5}) if and only if $\Om$ is a ball.

\vskip 0.1in

By Theorem 1.1 and a Faber-Krahn type inequality proved in section 3
Lemma 3.2, we have

\vskip 0.1in

{\bf Corollary 1.2.} Let $u(x)$ be the unique solution of problem
(\ref{4.1.6}) and $h(x)$ be the unique solution of problem
(\ref{4.1.7}). Then for any $k\geq q+1$, we have
\begin{eqnarray}\int_{\Om}u^k(x)dx\leq \int_{\Om^*}h^k(x)dx, \label{4.1.8}
\end{eqnarray}
 and
\begin{eqnarray}\max\limits_{x\in\Om}u(x)\leq \max\limits_{x\in\Om^*}h(x).\label{4.1.9}
\end{eqnarray}
Moreover, the equality holds in each of inequalities (\ref{4.1.8})
and (\ref{4.1.9}) if and only if $\Om$ is a ball.

\vskip 0.1in

Thanks to Corollary 1.2 and an explicit bound of solution of problem
(\ref{4.1.7}), we have

\vskip 0.1in

{\bf Corollary 1.3.} Let $u(x)$ be the unique solution of problem
(\ref{4.1.6}), and $\omega_n$ is the volume of unit ball in $R^n$.
Then
\begin{eqnarray}\max\limits_{x\in\Om}u(x)\leq \left[\frac{|\Om|}{\omega_n(2n)^{\frac{n}{2}}}
\right]^{\frac{2}{(1-q)n}}
\label{4.1.10}
\end{eqnarray}
with equality only if $\Omega$ is a ball.

 \vskip 0.1in

{\bf Remark 1.4.} Let $u(x)$ be the unique solution of problem
(\ref{4.1.6}). If $|\Om|<\omega_n(2n)^{\frac{n}{2}}$, then it
follows from Corollary 1.3 that $u(x)\rightarrow 0$ uniformly on
$\Omega$ when $q\rightarrow 1^-$. It is interest to know the
asymptotic behavior of $u(x)$ when
$|\Om|\geq\omega_n(2n)^{\frac{n}{2}}$ and $q\rightarrow 1^-$. It is
also interest to know the asymptotic behavior of $u(x)$ when
$q\rightarrow 0^+$.

\vskip 0.1in

{\bf Remark 1.5.}\ All results of this paper can be generalized to
p-Laplace equation with some modification of our method (see
\cite{ref DaiH}).

\vskip 0.1in

 The paper is organized as follows: As preliminary, we give some basic facts about the rearrangement of functions in section 2. In section 3, we prove a Chiti type
comparison result which is essential to the proof of our main
results. The proofs of Theorem 1.1, Corollary 1.2 and Corollary 1.3
are given in section 4.

 \vskip 1cm

 \section{Preliminary}

\setcounter{section}{2}

\setcounter{equation}{0}

\renewcommand{\theequation}{\thesection.\arabic{equation}}

\vskip 0.5cm

In this section, we recall some basic facts about the rearrangement
of functions and the existence and uniqueness result of problem
(\ref{4.1.6}).

Let $\Om$ be a bounded domain in $R^n$. The Schwartz symmetrization
$\Om^*$ of $\Om$ is a ball in $R^n$ with radius $R^*$ and centered
at $0$ such that $|\Om^*|=|\Om|$. Here, $|\Om|$ denotes the Lebesgue
measure of $\Om$. If we denote by $\om_n$ the volume of unit ball in
$R^n$, then it is easy to see
$$R^*=\left(\frac{|\Om|}{\om_n}\right)^{\frac{1}{n}}.$$

Let $f:\ \Om\mapsto R$ be a nonnegative measurable function. For any
$t\geq 0$. The level set $\Om_t$ of $f$ at the level $t$ is defined
by
$$\Om_t\doteq\{x\in \Om:\ f(x)>t\},\ \ \ \ t\geq 0.$$
The distribution function of $f$ is given by
$$\mu_f(t)=|\Om_t|=\mbox{meas}\{x\in \Om:\ f(x)>t\},\ \ \ \ t\geq 0.$$
Obviously, $\mu_f(t)$ is a monotonically decreasing function of $t$
and $\mu_f(t)=0$ for $t\geq \mbox{ess}.\sup.f$, while
$\mu_f(t)=|\Om|$ for $t=0$.

\vskip 0.1in

{\bf Definition 2.1.}  Let $\Om$ be a bounded domain in $R^n$, $f:\
\Om\mapsto R$ be a nonnegative measurable function. Then the
decreasing rearrangement $f^*$ of $f$ is a function defined on $[0,\
\infty)$ by
\begin{eqnarray*}
f^*(s)=\left
 \{\begin{array}{ll}
 \mbox{ess}.\sup.f &\mbox{for\ }s=0\\
 \\
 \inf\{t>0|\mu_f(t)<s\} &\mbox{for\ }s>0.
\end{array}
\right.
\end{eqnarray*}

Obviously, $f^*(s)=0$ for $s\geq |\Om|$. The increasing
rearrangement $f_*$ of $f$ is defined by $f_*(s)=f^*(|\Om|-s)$ for
$s\in(0,\ +\infty)$.

\vskip 0.1in

{\bf Definition 2.2.} Let $\Om$ be a bounded domain in $R^n$, $f:\
\Om\mapsto R$ be a nonnegative measurable function. Then the
decreasing Schwarz symmetrization $f^{\star}$ of $f$ is a function
defined by
$$f^{\star}(x)=f^{*}(\om_n|x|^n)\ \ \ \mbox{for\ }x\in\Om^*.$$

There are many fine properties of rearrangement. Here, we only
collect some important properties needed in this paper.

\vskip 0.1in {\bf Proposition 2.3.} Let $f:\ \Om\mapsto R$ be a
nonnegative measurable function. Then, $f,\ f^*$ and $f^{\star}$ are
all equimeasurable and
$$\int_{\Om}fdx=\int^{|\Om|}_0f^*(s)ds=\int_{\Om^*}f^{\star}(x)dx.$$
Moreover,  for any Borel measurable function $F:\ R\mapsto R$, there
holds

$$
\int_{\Om}F(f(x))dx=\int^{|\Om|}_0F(f^*(s))ds
=\int_{\Omega^*}F(f^{\star}(x))dx.$$

\vskip 0.1in {\bf Proposition 2.4.} If $f:\ [0,\ l]\mapsto R$ is
nonnegative and non-increasing, then $f=f^* \ a.e.$

\vskip 0.1in
{\bf Proposition 2.5.} If $\psi:\ R\mapsto R$ is a non-decreasing function, then
$$\psi(f^*)=(\psi(f))^*,\ \ \psi(f^{\star})=(\psi(f))^{\star}$$
for any nonnegative measurable function  $f:\ \Om\mapsto R$.

\vskip 0.1in
{\bf Proposition 2.6.} Let $f\in L^p(\Om),\ g\in L^q(\Om)$ with $\frac{1}{p}+\frac{1}{q}=1$. Then
$$
\int^{|\Om|}_0f^*(s)g_*(s)ds\leq\int_{\Om}f(x)g(x)dx\leq\int^{|\Om|}_0f^*(s)g^*(s)ds,
$$
$$
\int_{\Om^*}f^{\star}(x)g_{\star}(x)dx\leq\int_{\Om}f(x)g(x)dx\leq\int_{\Om^*}f^{\star}(x)g^{\star}(x)dx.
$$
Consequently
$$
\int_{E}f(x)dx\leq\int^{|E|}_0f^*(s)ds=\int_{E^*}f^{\star}(x)dx.
$$
for any measurable set $E\subset\Om$.

 \vskip 0.1in

{\bf Proposition 2.7.} If $f\in H^1_0(\Om)$, then $f^{\star}\in
H^1_0(\Om^*)$
 and
 $$\int_{\Om}|\nabla f|^2dx\geq
\int_{B_{R^*}(0)}|\nabla f^{\star}|^2dx$$ where $B_{R^*}(0)=\Om^*$.
Moreover, the equality holds if and only if $\Om$ is a ball.

\vskip 0.1in

The proof of all propositions mentioned above can be found in
\cite{ref SK, ref BK}.

\vskip 0.1in

{\bf Proposition 2.8 (\cite{ref GJG}).} Let $M,\ \alpha,\ \beta$ be
real numbers such that $0<\alpha\leq \beta$ and $M>0$. Let $f,\ g$
be real functions in $L^\beta([0,\ M])$. If the decreasing
rearrangements of $f$ and $g$ satisfy the inequality
$$\int^{s}_0f^{*^{\alpha}}(t)dt\leq\int^{s}_0g^{*{^{\alpha}}}(t)dt \ \ \ \mbox{for}\ \ \ s\in[0,\ M],$$
then
$$\int^{M}_0f^{*^{\beta}}(t)dt\leq\int^{M}_0g^{*^{\beta}}(t)dt.$$

\vskip 0.1in

The following result may be well known. However, for the reader's
convenience, we give a proof here.

\vskip 0.1in

{\bf Proposition 2.9.} Problem (\ref{4.1.6}) has an unique solution.

\vskip 0.1in

{\bf Proof.}\ Let $h(x)$ be the unique solution of
\begin{eqnarray*}
\left
 \{\begin{array}{ll}
  -\Delta h=1, &x\in \Om,\\
   h=0,& x\in\partial \Om.
\end{array}
\right.
\end{eqnarray*}
 Choose $M_0$ so that $M_0>M^q_0\max\limits_{x\in\Om}h^q(x)$, this is possible since $0<q<1$.

 Let $v_0(x)=M_0h(x)$, then
 $$-\Delta v_0=-M_0\Delta h(x)=M_0>M^q_0\max\limits_{x\in\Om}h^q(x)\geq v^q_0.$$
This implies that $v_0(x)$ is  sup-solution of problem (\ref{4.1.6}).

Let $\varphi_1(x)$ be the first eigenfunction of the eigenvalue problem
\begin{eqnarray*}
\left
 \{\begin{array}{ll}
  -\Delta\varphi=\lambda\varphi, &x\in \Om,\\
   \varphi=0,& x\in\partial \Om.
\end{array}
\right.
\end{eqnarray*}
We choose $\varphi_1(x)$ so that
$$
\varphi_1(x)>0,\ \ \max\limits_{x\in\Om}\varphi_1(x)=1.
$$

Let $v_{\eta_0}=\eta_0\varphi_1(x)$, then
$$
-\Delta v_{\eta_0}=-\eta_0\Delta\varphi_1(x)=\lambda_1\eta_0\varphi_1(x).
$$

Since $0<q<1$, we can choose $\eta_0$ small enough such that
$$
\lambda_1\eta_0\varphi_1(x)\leq\eta^q_0\varphi^q_1(x)=v^q_{\eta_0}.
$$
Hence $v_{\eta_0}$ is a sub-solution of problem (\ref{4.1.6}).
Choosing $\eta_0$ even more smaller, we can assume that
$v_{\eta_0}\leq v_0(x)$. Then by the sub- and super- solution
method, we know that problem (\ref{4.1.6}) has at least one solution
$u(x)$ which satisfies $v_{\eta_0}\leq u(x)\leq v_0(x)$.

To prove the uniqueness, we assume that $u_1(x)$ and $u_2(x)$ are
any two solutions of problem (\ref{4.1.6}). It is obvious that for
$b>0$ small enough, we have
$$
u_1(x)>bu_2(x),\ \ \ \ \ x\in\Om.
$$
Let
$$
b_0=\sup\{ b|\ u_1(x)>bu_2(x), \ \ x\in\Om\}.
$$
Then
$$
u_1(x)\geq b_0u_2(x),\ \ \ \ x\in\Om.
$$
and there exists at least one point $x_0\in\Om$ such that
\begin{eqnarray}
u_1(x_0)= b_0u_2(x_0),\ \ \ \ x\in\Om.\label{eq2.1}
\end{eqnarray}
If $b_0<1$, then $v_0=b_0u_2(x)$ satisfies
$$
-\Delta v_0=-b_0\Delta u_2(x)=b_0u^q_2(x)<b^q_0u^q_2(x)=v^q_0.
$$
Let $w=u_1(x)-v_0(x)$, then $w(x)$ satisfies
\begin{eqnarray*}
\left
 \{\begin{array}{ll}
  -\Delta w>u^q_1(x)-v^q_0(x)\geq0, &x\in \Om,\\
   w=0,& x\in\partial \Om.
\end{array}
\right.
\end{eqnarray*}
It follows from the strong maximum principle that $w(x)>0,\ x\in\Om$.
Hence
$$\ u_1(x)>b_0u_2(x),\ \ \ \ x\in\Om.$$
This contradicts (\ref{eq2.1}). Thus we must have $b_0\geq 1$ and
$$u_1(x)\geq u_2(x),\ \ \ \mbox{for}\ x\in\Om.$$
Changing the position of  $u_1(x)$ and $u_2(x)$, a similar argument implies that
$$u_2(x)\geq u_1(x),\ \ \ \mbox{for}\ x\in\Om.$$
Consequently,
$$u_1(x)\equiv u_2(x),\ \ \ \mbox{for}\ x\in\Om.$$
This means that problem (\ref{4.1.6}) has only one solution.

\vskip 1cm

\section{Chiti Type Comparison Result}

\setcounter{section}{3}

\setcounter{equation}{0}

\renewcommand{\theequation}{\thesection.\arabic{equation}}

\vskip 0.5cm

 Let $\Om$ be a bounded domain in $R^n$, and $\|\cdot\|_{L^{q+1}(\Om)}$ denote the norm of space $L^{q+1}(\Om)$. We define
 $$
 S_q(\Om)=\inf\limits_{v\in H^1_0(\Om)}\{\int_{\Om}|\nabla v|^2dx\biggr|\ \|v\|^2_{L^{q+1}(\Om)}=1. \}
 $$

It is easy to prove that $ S_q(\Om)$ can be achieved by an unique
positive function $v(x)$. Moreover, $v(x)$ satisfies
  \begin{eqnarray}
\left
 \{\begin{array}{ll}
-\Delta v(x)=S_q(\Om)v^q(x), &x\in \Om,\\
v(x)>0,& x\in\Om,\\
v(x)=0,& x\in\partial \Om,\\
\int_{\Om}v^{q+1}(x)dx=1.
\end{array}
\right.\label{eq3.1}
\end{eqnarray}

In this section, we prove a Chiti type comparison result for problem
(\ref{eq3.1}). To this end, we need some lemmas first.

\vskip 0.1in

{\bf Lemma 3.1.} For any $\lambda>0$ and $\lambda\neq S_q(\Om)$, the
following problem has no solution
  \begin{eqnarray}
\left
 \{\begin{array}{ll}
-\Delta f(x)=\lambda f^q(x), &x\in \Om,\\
f(x)>0,& x\in\Om,\\
f(x)=0,& x\in\partial \Om,\\
\int_{\Om}f^{q+1}(x)dx=1.
\end{array}
\right.\label{eq3.2}
\end{eqnarray}

{\bf Proof.} We prove Lemma 3.1 by contradiction. Assume that
problem (\ref{eq3.2}) has a solution $f_{\lambda_0}$ for some
$\la_0>0$ and $\la_0\neq S_q(\Om)$. Then, it is easy to check that
$\widetilde f=\la^{\frac{1}{q-1}}_0f_{\la_0}$ is a solution of
problem (\ref{4.1.6}) which satisfies
$$
\int_{\Om}{\widetilde f}^{q+1}(x)dx=\la^{\frac{q+1}{q-1}}_0.
$$
On the other hand, if we denote by $v(x)$ the minimizer of $S_q(\Om)$, then $\widetilde v=S^{\frac{1}{q-1}}_q(\Om)v(x)$ is
also a solution of problem (\ref{4.1.6}) which satisfies
$$
\int_{\Om}{\widetilde v}^{q+1}(x)dx=S^{\frac{q+1}{q-1}}_q(\Om).
$$
It is obvious that $\widetilde v\neq \widetilde f$ due to $\la_0\neq
S_q(\Om)$. Hence problem (\ref{4.1.6}) has at least two solutions
$\widetilde v$ and $\widetilde f$. This contradicts Proposition 2.9.

\vskip 0.1in

{\bf Lemma 3.2.} $S_q(\Om)\geq S_q(\Om^*)$ with equality if and only
if $\Om$ is a ball.

\vskip 0.1in

{\bf Proof.} Let $v(x)$ be the minimizer of $S_q(\Om)$ and
$v^{\star}(x)$ be its Schwartz symmetrization. Then by Proposition
2.3 and Proposition 2.7, we have
 $$\int_{\Om}|\nabla v|^2dx\geq
\int_{\Om^*}|\nabla v^{\star}|^2dx,$$
$$
\int_{\Om}v^{q+1}(x)dx=\int_{\Om^*}{(v^{\star})}^{q+1}(x)dx=1.
$$
Hence, by the definition of $S_q(\Om^*)$, we have
$$S_q(\Om^*)\leq\int_{\Om^*}|\nabla v^{\star}|^2dx\leq\int_{\Om}|\nabla v|^2dx=S_q(\Om).$$
If $S_q(\Om^*)=S_q(\Om)$, then $\int_{\Om^*}|\nabla
v^{\star}|^2dx=\int_{\Om}|\nabla v|^2dx$. Hence, by Proposition 2.7,
we know that $\Om$ is a ball.

\vskip 0.1in

Let $\sigma_3=\frac{q+1}{n+2-(n-2)q}$. Then the following lemma
holds

\vskip 0.1in

{\bf Lemma 3.3.} Let $v(x)$ be the minimizer of $S_q(\Om^*)$ and $r_{*}=\left(\frac{S_q(\Om^*)}{S_q(\Om)}\right)
^{\sigma_3}R^*$. Then $S_q(B_{r_*}(0))=S_q(\Om)$ and the minimizer of $S_q(B_{r_*}(0))$ is $z(y)=
\left(\frac{R^*}{r_*}\right)^{\frac{n}{q+1}}v(\frac{R^*}{r_*}y)$ for $y\in B_{r_*}(0)$.

\vskip 0.1in

{\bf Proof.} Since $v(x)$ is the minimizer of $S_q(\Om^*),\ v(x)$
satisfies
\begin{eqnarray*}
\left\{\begin{array}{ll}
-\Delta v(x)=S_q(\Om^*)v^q(x), &x\in \Om^*,\\
v(x)>0,& x\in\Om^*,\\
v(x)=0,& x\in\partial \Om^*,\\
\int_{\Om^*}v^{q+1}(x)dx=1.
\end{array}
\right.
\end{eqnarray*}
Let $x=\frac{R^*}{r_*}y$ and $H(y)=v(\frac{R^*}{r_*}y)$. Then
 \begin{eqnarray*}
\frac{\partial H}{\partial y_i}=\frac{R^*}{r_*}\frac{\partial v}{\partial x_i},\\
\frac{\partial^2 H}{\partial y^2_i}=(\frac{R^*}{r_*})^2\frac{\partial^2 v}{\partial x^2_i}.
\end{eqnarray*}
Hence
 \begin{eqnarray*}
 -\Delta H(y)=-(\frac{R^*}{r_*})^2\Delta v=(\frac{R^*}{r_*})^2S_q(\Om^*)H^q(y),\ \ y\in B_{r_*}(0).
\end{eqnarray*}
Noting that
\begin{eqnarray*}
1=\int_{\Om^*}v^{q+1}(x)dx&=&(\frac{R^*}{r_*})^n\int_{B_{r_*}(0)}H^{q+1}(y)dy\\
&=&\int_{B_{r_*}(0)}\left[(\frac{R^*}{r_*})^{\frac{n}{q+1}}H(y)\right]^{q+1}dy,
\end{eqnarray*}
if we let
$z(y)=(\frac{R^*}{r_*})^{\frac{n}{q+1}}H(y)=(\frac{R^*}{r_*})^{\frac{n}{q+1}}v(\frac{R^*}{r_*}y)$,
then $z(y)$ satisfies
  \begin{eqnarray*}
\left
 \{\begin{array}{ll}
-\Delta z(y)=\left(\frac{R^*}{r_*}\right)^{\frac{1}{\sigma_3}}S_q(\Om^*)z^q(y), &y\in B_{r_*}(0),\\
z(y)>0,& y\in B_{r_*}(0),\\
z(y)=0,& y\in\partial B_{r_*}(0),\\
\int_{B_{r_*}(0)}z^{q+1}(y)dy=1.
\end{array}
\right.
\end{eqnarray*}
Hence, by Lemma 3.1, we have
\begin{eqnarray*}
S_q(B_{r_*}(0))=\left(\frac{R^*}{r_*}\right)^{\frac{1}{\sigma_3}}S_q(\Om^*)=S_q(\Om).
\end{eqnarray*}
and the minimizer of $S_q(B_{r_*}(0))$ is $z(y)=
\left(\frac{R^*}{r_*}\right)^{\frac{n}{q+1}}v(\frac{R^*}{r_*}y)$.
This completes the proof of Lemma 3.3.

\vskip 0.1in

By Lemma 3.2 and the definition of $r_*$, we have
$B_{r_*}(0)\subset\Om^*$ with equality if and only if $\Om$ is a
ball. Let $M=|\Om|$ and $M_*=|B_{r_*}(0)|$, then $M_*\leq M$. The
main result of this section is the following Chiti type comparison
result.

\vskip 0.1in {\bf Theorem 3.4.} Let $v(x)$ be the minimizer of
$S_q(\Om)$ and $z(x)$ be the minimizer of $S_q(B_{r_*}(0))$. If we
denote by $v^*(s)$ the decreasing rearrangement of $v(x)$, and
$z^*(s)$ the decreasing rearrangement of $z(x)$, then there exists
an unique point $s_0\in (0,\ M_*)$ such that
  \begin{eqnarray*}
\left
 \{\begin{array}{ll}
z^*(s)>u^*(s) &\mbox{for}\ s\in [0,\ s_0)\\
\\
z^*(s)<u^*(s) &\mbox{for}\ s\in (s_0,\ M_*].
\end{array}
\right.
\end{eqnarray*}

{\bf Proof.} Since $u(x)$ is the minimizer of $S_q(\Om)$, it is easy to see that $u(x)$
 satisfies
  \begin{eqnarray}
\left
 \{\begin{array}{ll}
-\Delta u(x)=S_q(\Om)u^q(y), &x\in \Om,\\
u(x)>0,&x\in \Om,\\
u(x)=0,& x\in\partial\Om.\\
\end{array}
\right.\label{eq3.3}
\end{eqnarray}
From this, we can prove that the decreasing rearrangement $u^*(s)$ of $u(x)$ satisfies
\begin{eqnarray}
-\frac{du^*(s)}{ds}\leq S_q(\Om)n^{-2}\omega^{\frac{-2}{n}}_n
s^{-\frac{2(n-1)}{n}}\int^s_0{(u^*)^q(t)}dt\ \ \ \  a.e.\ \mbox{in}\
[0,\ M],\label{eq3.4}
\end{eqnarray}
In fact, integrating the first equation in (\ref{eq3.3}) over
$\Om_t=\{x\in\Om\ |\ u(x)>t\}$, we have
\begin{eqnarray}
  -\int_{\partial\Om_t}\frac{\partial u(x)}{\partial\nu}ds=S_q(\Om)\int_{\Om_t}u^qdx.\label{eq3.5}
\end{eqnarray}
Since $\partial\Om_t=\{x\in\\Omega \ |\ u(x)=t\}$, we have
\begin{eqnarray}
  -\int_{\partial\Om_t}\frac{\partial u(x)}{\partial\nu}ds=\int_{\partial\Om_t}|\nabla u|ds.\label{eq3.6}
\end{eqnarray}
Noting that
\begin{eqnarray*}
\int_{\partial\Om_t}|\nabla u|ds\int_{\partial\Om_t}\frac{ds}{|\nabla u|}\geq|\partial\Om_t|^2.
\end{eqnarray*}
It follows from the isoperimetric inequality
\begin{eqnarray}
 \int_{\partial\Om_t}|\nabla u|ds\int_{\partial\Om_t}\frac{ds}
 {|\nabla u|}\geq n^2\omega^{\frac{2}{n}}_n|\Om_t|^{\frac{2(n-1)}{n}}.\label{eq3.7}
\end{eqnarray}
By Co-area formula, we have
\begin{eqnarray*}\mu(t)=|\Om_t|=\int_{\Om_t}dx=\int^{+\infty}_t\int_{\partial\Om_t}\frac{ds}
 {|\nabla u|}.
 \end{eqnarray*}
Consequently,
\begin{eqnarray}
\frac{d\mu(t)}{dt}=-\int_{\partial\Om_t}\frac{ds}{|\nabla u|}.\label{eq3.8}
\end{eqnarray}
From (\ref{eq3.5}), (\ref{eq3.6}), (\ref{eq3.7}) and (\ref{eq3.8}), we obtain
\begin{eqnarray}
\frac{n^2\omega^{\frac{2}{n}}_n(\mu(t))^{\frac{2(n-1)}{n}}}{-\mu^{'}(t)}\leq S_q(\Om)\int_{\Om_t}u^qdx.
\label{eq3.9}
\end{eqnarray}
Since $\Om_t\subset\Om$, we have
\begin{eqnarray}
\int_{\Om_t}u^qdx\leq\int^{|\Om_t|}_{0}(u^q)^*(\tau)d\tau=\int^{\mu(t)}_{0}(u^*(\tau))^qd\tau.
\label{eq3.10}
 \end{eqnarray}
Combing (\ref{eq3.9}) with (\ref{eq3.10}), we obtain
\begin{eqnarray*}
-\frac{1}{\mu^{'}(t)}\leq S_q(\Om)n^{-2}\omega^{\frac{-2}{n}}_n
(\mu(t))^{-\frac{2(n-1)}{n}}\int^{\mu(t)}_0(u^*(\tau))^{q}d\tau.
\end{eqnarray*}
Noticing that $u^*(s)$ is essentially an inverse of $\mu(t)$, we
have
 \begin{eqnarray*}
 -\frac{du^*(s)}{ds}\leq S_q(\Om)n^{-2}\omega^{\frac{-2}{n}}_n
s^{-\frac{2(n-1)}{n}}\int^{s}_0(u^*)^q(\tau))d\tau.
 \end{eqnarray*}
This is just the desired conclusion of (\ref{eq3.4}).

Since $S_q(B_{r_*}(0))=S_q(\Om)$, the minimizer $z(x)$ of
$S_q(B_{r_*}(0))$ satisfies
\begin{eqnarray}
\left
 \{\begin{array}{ll}
-\Delta z(x)=S_q(\Om)z^q(x), &x\in B_{r_*}(0),\\
z(x)>0,& x\in B_{r_*}(0),\\
z(x)=0,& x\in\partial B_{r_*}(0).
\end{array}
\right.\label{eq3.11}
\end{eqnarray}
Noticing that uniqueness result valid for (\ref{eq3.11}), it is
trivial to see that $z$ is radial symmetry. That is $z(x)=z(|x|)$.
Moreover, as a function of $s=\omega_n|x|^n$, $z(s)$ is decreasing.
Hence, by making use of (\ref{eq3.11}), Proposition 2.4 and
Proposition 2.5, a similar argument to that used to derive
(\ref{eq3.4}) implies that
\begin{eqnarray}
-\frac{dz^*(s)}{ds}= S_q(\Om)n^{-2}\omega^{\frac{-2}{n}}_n
s^{-\frac{2(n-1)}{n}}\int^s_0(z^*)^q(t)dt\ \ \ \  a.e.\ \mbox{in}\
[0,\ M_*],\label{eq3.12}
\end{eqnarray}
Now, Theorem 3.3 can be proved by making use of (\ref{eq3.4}) and
(\ref{eq3.12}). To this end, we first note that there exists at
least one point $s_0\in (0,\ M_*)$ such that $u^*(s_0)=z^*(s_0)$
because of
\begin{eqnarray*}
\int_{\Om}u^{q+1}(x)dx=\int^{M_*}_0(u^*)^{q+1}(s)ds=1=\int^{M_*}_0(z^*(s))^{q+1}ds=\int_{B_{r_*}(0)}z^{q+1}(x)dx.
 \end{eqnarray*}
Next, we prove that there exists only one point $s_0\in (0,\ M_*)$ such that $u^*(s_0)=z^*(s_0)$. Otherwise, there would exist at
least two points $s_1,\ s_2\in (0,\ M_*)$ such that
\begin{eqnarray*}
u^*(s_1)=z^*(s_1),\ \ \ u^*(s_2)=z^*(s_2).
\end{eqnarray*}
This would imply that there exists an interval $[s_1,\ s_2]\subset [0,\ M_*)$ such that
\begin{eqnarray*}
\left
\{\begin{array}{lll}
u^*(s_i)=z^*(s_i),&i=1,\ 2;\\
u^*(s)>z^*(s), &s\in(s_1,\ s_2).
 \end{array}
\right.
 \end{eqnarray*}
 Let
\begin{eqnarray*}
w(s)= \left
\{\begin{array}{lll}
z^*(s),  & \mbox{if}\ \int^s_0(u^*(\tau))^{q}d\tau\leq\int^{s}_0(z^*(\tau))^{q}d\tau,\ s\in[0,\ s_1];\\
u^*(s),  & \mbox{if}\ \int^s_0(u^*(\tau))^{q}d\tau\geq\int^{s}_0(z^*(\tau))^{q}d\tau,\ s\in[0,\ s_1];\\
u^*(s), & s\in[s_1,\ s_2];\\
z^*(s), & s\in[s_2,\ M_*].
 \end{array}
\right.
\end{eqnarray*}
Then, it is easy to verify that $w(s)$ satisfies
\begin{eqnarray}
\left\{\begin{array}{lll}
-\frac{d w(s)}{ds}\leq S_q(\Om)n^{-2}\omega_n^{-\frac{2}{n}}
s^{-\frac{2(n-1)}{n}}\int^s_0 w^q(t)dt,& a.e.\ \mbox{in}\ [0,\ M_*],\\
w(s)>0,&s\in (0,\ M_*),\\
w(M_*)=0,\\
\|w\|_{L^{q+1}(0,\ M_*)}\geq 1.
\end{array}
\right.\label{eq3.13}
\end{eqnarray}
Define $$\eta(x)=
\frac{w(\omega_n|x|^n)}{\|w(\omega_n|x|^n)\|_{q+1(B_{r_*}(0))}}.$$
Then, $\eta(x)\in W^{1,\ 2}_0(B_{r_*}(0))$ and
$\|\eta(x)\|_{q+1(B_{r_*}(0))}=1$. Since $\eta(x)$ is obviously not
the minimizer of $S_q(B_{r_*}(0))$, we have
$$S_q(\Om)=S_q(B_{r_*}(0))<\int_{B_{r_*}(0)}|\nabla\eta(x)|^2dx.$$
Since
\begin{eqnarray*}
\int_{B_{r_*}(0)}|\nabla\eta(x)|^2dx&=&n^{2}\omega^{\frac{2}{n}}_n\int^{M_*}_0|\eta^{'}(s)|^2s^{\frac{2(n-1)}{n}}ds\\
&=&\frac{n^{2}\omega^{\frac{2}{n}}_n}{\|w\|^2_{q+1(B_{r_*}(0))}}\int^{M_*}_0|w^{'}(s)|^2s^{\frac{2(n-1)}{n}}ds,
\end{eqnarray*}
and
\begin{eqnarray*}
n^{2}\omega^{\frac{2}{n}}_n\int^{M_*}_0|w^{'}(s)|^2s^{\frac{2(n-1)}{n}}ds&=&
n^{2}\omega^{\frac{2}{n}}_n\int^{M_*}_0(-w^{'}(s))(-w^{'}(s))s^{\frac{2(n-1)}{n}}ds\\
&\leq&S_q(\Om)\int^{M_*}_0(-w^{'}(s))\int^s_0w^q(\tau)d\tau ds\\
&=&S_q(\Om)\int^{M_*}_0w^{q+1}(s)ds\\
&=&S_q(B_{r_*}(0))\|w\|^{q+1}_{q+1(B_{r_*}(0))}
\end{eqnarray*}
We have
$$
\int_{B_{r_*}(0)}|\nabla\eta(x)|^2dx\leq S_q(B_{r_*}(0))\|w\|^{q+1-2}_{q+1(B_{r_*}(0))}=S_q(B_{r_*}(0))\|w\|^{q-1}_{q+1(B_{r_*}(0))}.
$$
Thus
$$S_q(B_{r_*}(0))<\int_{B_{r_*}(0)}|\nabla\eta(x)|^2dx\leq S_q(B_{r_*}(0))\|w\|^{q-1}_{q+1(B_{r_*}(0))}.$$
Noticing that $\|w\|_{q+1(B_{r_*}(0))}\geq 1$ and $q-1<0$, we obtain
$$S_q(B_{r_*}(0))<S_q(B_{r_*}(0)).$$
This is a contradiction.

Hence, there exists only one point $s_0\in (0,\ M_*)$ such that $z^*(s_0)=u^*(s_0)$ and this implies that
\begin{eqnarray*}
\left
\{\begin{array}{lll}
z^*(s)>u^*(s),&\mbox{for}\ s\in (0,\ s_0)),\\
\\
z^*(s)<u^*(s), &\mbox{for}\ s\in (s_0,\ M_*).
\end{array}
\right.
 \end{eqnarray*}
 So, we complete the proof of Theorem 3.3.

\vskip 0.1in

{\bf  Corollary 3.4.} Let $u(x)$ be the minimizer of $S_q(\Om)$ and
$z(x)$ be the minimizer of $S_q(B_{r_*}(0))$. Then for any $k\geq
q+1$, there holds
$$\int_{\Om}u^kdx\leq \int_{B_{r_*}(0)}z^k(x)dx.$$
It follows that
$$\sup\limits_{x\in\Om}u(x)\leq\sup\limits_{x\in B_{r_*}(0)}z(x).$$
Moreover, the equality holds in the above two inequalities if and
only if $\Omega$ is a ball.

\vskip 0.1in

{\bf Proof.}\ By the proposition of rearrangement, we have
 $$
 \int^{M}_0(u^*)^{q+1}(s)ds=1=\int^{M_*}_0(z^*(s))^{q+1}ds.
 $$
 Hence
$$
 \int^{M_*}_0(u^*)^{q+1}(s)ds\leq\int^{M_*}_0(z^*(s))^{q+1}ds.
 $$
 Let $s_0$ be the point in $(0,\ M_*)$ determined in Theorem 3.3.
 Then
 $$
 \int^{M_*}_{s_0}(u^*)^{q+1}(s)ds-\int^{M_*}_{s_0}(z^*)^{q+1}(s)ds\leq \int^{s_0}_0\left((z^*)^{q+1}-(u^*)^{q+1}\right)(s)ds.
 $$
 Since $u^*(s)\geq z^*(s)$ for any $s\in[s_0,\ M_*]$. It follows that for any $s\in[s_0,\ M_*]$, there holds
 $$
 \int^s_{s_0}\left((u^*)^{q+1}-(z^*)^{q+1}\right)(s)ds\leq\int^{s_0}_0
 \left((z^*)^{q+1}-(u^*)^{q+1}\right)(s)ds
 $$
Consequently,
$$
 \int^s_0(u^*)^{q+1}(\tau)d\tau\leq\int^s_0(z^*)^{q+1}(\tau)d\tau\ \ \ \mbox{for\ any\ }s\in (0,\ M_*).
$$
By the definition of $z^*(s)$, we have $z^*(s)=0$ for $s\geq M_*)$.
Hence
$$
 \int^s_0(u^*)^{q+1}(\tau)d\tau\leq\int^s_0(z^*)^{q+1}(\tau)d\tau\ \ \ \mbox{for\ any\ }s\in (0,\ M).
$$
From this and Proposition 2.8, we have
$$
\int^{M}_0(u^*)^k(s)ds\leq\int^{M_*}_0(z^*)^k(s)ds.
$$
Noticing that
$$
\int_{\Om}u^k(x)dx=\int^{M}_0(u^*)^k(s)ds,
$$
$$
\int_{B_{r_*}(0)}z^k(x)dx=\int^{M_*}_0(z^*)^k(s)ds.
$$
We obtain
$$
\int_{\Om}u^k(x)dx\leq\int_{B_{r_*}(0))}z^k(x)dx
$$
for any $k\geq q+1$. This completes the proof of Corollary 3.4.
%\begin{eqnarray}  \end{eqnarray}         \begin{eqnarray*}  \end{eqnarray*}

\vskip 1cm

\section{Proofs of Theorem 1.1, Corollary 1.2 and Corollary 1.3}

\setcounter{section}{4}

\setcounter{equation}{0}

\renewcommand{\theequation}{\thesection.\arabic{equation}}

\vskip 0.5cm

In this section, we prove Theorem 1.1, Corollary 1.2 and Corollary
1.3. For simplicity, we always use the notations $\sigma_1,
\sigma_2$ and $\sigma_3$ introduced in section 1 and section 3 in
this section.

\vskip 0.1in

{\bf Proof of Theorem 1.1.}\ Let $u(x)$ be the solution of problem
(\ref{4.1.6}). Then $v(x)=\frac{u(x)}{\|u\|_{L^{q+1}(\Om)}}$
satisfies
  \begin{eqnarray*}
\left
 \{\begin{array}{ll}
-\Delta v(x)=\|u\|^{q-1}_{L^{q+1}(\Om)}v^q(x), &x\in \Om,\\
v(x)>0,& x\in\Om,\\
v(x)=0,& x\in\partial \Om,\\
\int_{\Om}v^{q+1}(x)dx=1.
\end{array}
\right.
\end{eqnarray*}
Hence, by Lemma 3.1, we have $S_q(\Om)=\|u\|^{q-1}_{L^{q+1}(\Om)}$ and the minimizer of $S_q(\Om)$ is $v(x)$.

Similarly, if $h(x)$ is the unique solution of problem (\ref{5}). Then $S_q(\Om^*)=\|h\|^{q-1}_{L^{q+1}(\Om^*)}$ and
the minimizer of $S_q(\Om^*)$ is $\frac{h(x)}{\|h\|_{L^{q+1}(\Om^*)}}$.

By the definition of $r_*$, we have
$r_*=\left[\frac{\|h\|_{L^{q+1}(\Om^*)}}{\|u\|_{L^{q+1}(\Om)}}\right]^
{(q-1)\sigma_3}R^*$. Moreover, by Lemma 3.3, we know that the
minimizer of $S_q(B_{r_*}(0))$ is
$$z(x)=\left(\frac{R^*}{r_*}\right)^{\frac{n}{q+1}}\frac{h(\frac{R^*}{r_*}x)}{\|h\|_{L^{q+1}(\Om^*)}}.$$
Applying Corollary 3.4 to $v(x)$ and $z(x)$, we have, for any $k\geq
q+1$, that
\begin{eqnarray*}
\int_{\Om}u^k(x)dx&\leq&\frac{\|u\|^k_{L^{q+1}(\Om)}}{\|h\|^k_{L^{q+1}(\Om^*)}}\int_{B_{r_*}(0)}
\left(\frac{R^*}{r_*}\right)^{\frac{nk}{q+1}}h^k(\frac{R^*}{r_*}x)dx\\
&=&\frac{\|u\|^k_{L^{q+1}(\Om)}}{\|h\|^k_{L^{q+1}(\Om^*)}}\left(\frac{R^*}{r_*}\right)^{\frac{nk}{q+1}-n}\int_{\Om^*}h^k(x)dx
\end{eqnarray*}
Since
$$
\left(\frac{R^*}{r_*}\right)^{\frac{nk}{q+1}-n}=\left(\frac{R^*}{r_*}\right)^{\frac{(k-q-1)n}{q+1}}=\left[\frac{\|h\|_{L^{q+1}(\Om^*)}}{\|u\|_{L^{q+1}(\Om)}}\right]^
{-\frac{(k-q-1)(q-1)n}{n+2-(n-2)q}}
$$
We have
\begin{eqnarray*}
\int_{\Om}u^k(x)dx&\leq&\frac{\|u\|^k_{L^{q+1}(\Om)}}{\|h\|^k_{L^{q+1}(\Om^*)}}
\left[\frac{\|h\|_{L^{q+1}(\Om^*)}}{\|u\|_{L^{q+1}(\Om)}}\right]^
{\frac{(k-q-1)(1-q)n}{n+2-(n-2)q}}\int_{\Om^*}h^k(x)dx.
\label{4.4.1}
\end{eqnarray*}
If we set
$$
C(q,k,\Om^*)=\int_{\Om^*}h^k(x)dx\biggr/\|h\|^{\sigma_1}_{L^{q+1}(\Om^*)},
$$
then
\begin{eqnarray*}
\int_{\Om}u^k(x)dx&\leq&C(q,k, \Om^*)\|u\|^{\sigma_1}_{L^{q+1}(\Om)}\\
\end{eqnarray*}
and the equality holds if and only $\Om$ is a ball.

\vskip 0.05in

If we set
$$C(q, \Om^*)=\mbox{ess.}\sup\limits_{x\in\Om^*}h(x)\biggr/\|h\|^{\sigma_2}_{L^{q+1}(\Om^*)},$$
then we can obtain
$$
\mbox{ess.}\sup\limits_{x\in\Om}u(x)\leq C(q,
\Om^*)\|u\|^{\sigma_2}_{L^{q+1}(\Om)}.
$$
and the equality holds if and only if $\Omega$ is a ball. This
completes the proof of Theorem 1.1.

\vskip 0.1in

{\bf Proof of Corollary 1.2.} Following the argument of theorem 1.1,
we know that for any $k\geq q+1$,
\begin{eqnarray}
\int_{\Om}u^k(x)dx&\leq&\frac{\|u\|^k_{L^{q+1}(\Om)}}{\|h\|^k_{L^{q+1}(\Om^*)}}
\left[\frac{\|h\|_{L^{q+1}(\Om^*)}}{\|u\|_{L^{q+1}(\Om)}}\right]^
{\frac{(k-q-1)(1-q)n}{n+2-(n-2)q}}\int_{\Om^*}h^k(x)dx.
\label{4.4.1}
\end{eqnarray}
Since $S_q(\Om)=\|u\|^{q-1}_{L^{q+1}(\Om)}$ and
$S_q(\Om^*)=\|h\|^{q-1}_{L^{q+1}(\Om^*)}$, we have
\begin{eqnarray}
\int_{\Om}u^k(x)dx&\leq&\frac{S^{\frac{k}{q-1}}_q(\Om)}{S^{\frac{k}{q-1}}_q(\Om^*)}
\left[\frac{S^{\frac{1}{q-1}}_q(\Om^*)}{S^{\frac{1}{q-1}}_q(\Om)}\right]^
{\frac{(k-q-1)(1-q)n}{n+2-(n-2)q}}\int_{\Om^*}h^k(x)dx.
\label{4.4.1}
\end{eqnarray}
Noting that $0<q<1$, it follows from Lemma 3.2 that
\begin{eqnarray*}
\int_{\Om}u^k(x)dx&\leq&\int_{\Om^*}h^k(x)dx.
\end{eqnarray*}
Consequently
\begin{eqnarray*}\max\limits_{x\in\Om}u(x)\leq \max\limits_{x\in\Om^*}h(x).
\end{eqnarray*}
Moreover, the equality in the above two inequalities holds if and
only if $\Om$ itself is a ball. This completes the proof of
Corollary 1.2.

\vskip 0.1in

{\bf Proof of Corollary 1.3.} Let $G(x,y)$ denote the Green's
function related to the Laplace operator on $\Om^*$. Then Green's
formula implies that the solution of problem (\ref{4.1.7}) can be
represented as
$$h(x)=\int_{\Om^*}G(x,y)h^q(y)dy.$$
Hence
\begin{eqnarray*}
\max\limits_{x\in\Om^*}h(x)\leq \max\limits_{x\in\Om^*}h^q(x)
\times\max\limits_{x\in\Om^*}\int_{\Om^*}G(x,y)dy
\leq\left[\max\limits_{x\in\Om^*}h(x)\right]^q
\times\max\limits_{x\in\Om^*}\int_{\Om^*}G(x,y)dy.
\end{eqnarray*}
Consequently
\begin{eqnarray*}
\max\limits_{x\in\Om}h(x)\leq\left[\max\limits_{x\in\Om^*}\int_{\Om^*}G(x,y)dy\right]
^{\frac{1}{1-q}}.
\end{eqnarray*}

Let $\gamma(x)=\int_{\Om^*}G(x,y)dy$. Then it is easy to verify that
$\gamma(x)$ satisfies
\begin{eqnarray}
\left \{\begin{array}{ll}
-\Delta \gamma=1 & \mbox{in}\ \Om^*\\
\gamma=0 & \mbox{on}\ \partial \Om^*.
\end{array}
\right.
 \label{4.4.3}
\end{eqnarray}
Hence, an easy computation tells that
$\gamma(x)=\left[\frac{|\Om^*|}{\omega_n(2n)^{\frac{n}{2}}}
\right]^{\frac{2}{n}}-\frac{1}{2n}|x|^2$ and
$\max\limits_{x\in\Om^*}\gamma(x)=
\left(\frac{|\Om|}{\omega_n(2n)^{\frac{n}{2}}}
\right)^{\frac{2}{n}}$. Thus we have
\begin{eqnarray*}
\max\limits_{x\in\Om}h(x)\leq
\left[\max\limits_{x\in\Om^*}\gamma(x)\right]^{\frac{1}{1-q}}
=\left[\frac{|\Om|}{\omega_n(2n)^{\frac{n}{2}}}
\right]^{\frac{2}{(1-q)n}}.
\end{eqnarray*}
Now, the conclusion of Corollary 1.3 follows from Corollary 1.2.
This completes the proof of Corollary 1.3.

% Set the ending of a LaTeX document

\begin{thebibliography}{CC}

\bibitem{ref ATL1}A. Alvino, G. Trombetti and P.-L. Lions,
Comparasion results for elliptic and parabolic equations via Schwarz
symmetrization, Ann. Inst. Henri Poincar\'{e}, Analyse non
lin\'{e}aire, 7(1990), p37-65.

\bibitem{ref ATL2}A. Alvino, P.-L. Lions and G. Trombetti,
Comparasion results for elliptic and parabolic equations via Schwarz
symmetrization: A new approach, Differ. Integral Equations, 4(1992),
p25-50.

\bibitem{ref MSA}M. S. Ashbaugh, R. D. Benguria, A sharp bound for the ratio of the first two
eigenvalue of Dirichlet Laplacians and extensions, Ann. Math.,
135(1992), p601-628.

\bibitem{ref CB}C. Bandle, Isoperimetric inequality for some
             eigenvalues of an inhomogeneous free membrane, SIAM
             J. Appl. Math., 22(1972), p142-147.

\bibitem{ref Ban1}C. Bandle, On symmetrization in parabolic equations, J. Anal. Math., 30(1976), p98-112.

\bibitem{ref MVB}M. Belloni, V. Ferone and B. Kawohl, Isoperimetric inequalities, Wulff shape
              and related questions for strongly nonlinear elliptic operators,
              Z. angew. Math. Phys., 54(2003) p771-783.

\bibitem{ref Bha} T. Bhattacharya, A proof of the Faber-Krahn inequality
for the first eigenvalue of the p-Laplacian, Annali di Matematica
pura ed applicata, IV 177(1999), p325-343.

\bibitem{ref Bos} M. H. Bossel, Membranes $\acute{e}$lastiquement
li$\acute{e}$es inhomog$\grave{e}$nes ou sur une surface: une
nouvelle extension du th$\acute{e}$or$\acute{e}$me
isop$\acute{e}$rim$\acute{e}$trique de Rayleigh-Faber-Krahn, Z.
Angew. Math. Phys., 39(1988), p733-742.

\bibitem{ref FC}F. Chiacchio, Estimates for the first eigenfunction
             of linear eigenvalue problem via steiner symmetrization, Publ. Mat.
             53(2009), p47-71.

\bibitem{ref GC}G. Chiti, A reverse H$\ddot{o}$lder inequality
                for the eigenfunctions of linear second order
                elliptic operators, Journal of Applied Mathematics and
                Physics, 33(1982), p143-148.


\bibitem{ref QD}Qiuyi Dai and Yuxia Fu, Faber-Krahn inequality for Robin problem involving p-
                Laplacian, 2008 Preprint, see arXiv:0912.0393.

\bibitem{ref DaiH}Qiuyi Dai and Huaxiang Hu, Isoperimetric
inequality and sharp bound for positive solution of p-Laplace
equation, In preparation.

\bibitem{ref Dan} D. Daners, A Faber-Krahn inequality for Robin problems
in any space dimension, Math. Ann. 335(2006), p767-785.

\bibitem{ref Fab} G. Faber, Beweis, dass unter allen homogenen Membranen
von gleicher Fl\"{a}che und gleicher Spannung die kreisf\"{o}rmige
den tiefsten Grundton gibt, Sitzungsber, Bayr. Akad. Wiss.
M\"{u}nchen, Math.-Phys. Kl(1923), p169-127.

\bibitem{ref GJG} G. H. Hardy, J. E. Littlewood and G. P$\acute{o}$lya,
Some simple inequalities satisfied by convex functions. Messenger
Math., 58(1929), p152.

\bibitem{ref BV}B. Kawohl, V. Fridman, Isoperimetric estimates for the first eigenvalue
               of the p-Laplace operator and the Cheeger constant, Comment.
               Math. Univ. Carolinae, 44(2003), p659-667.

\bibitem{ref BK}B. Kawohl, Rearrangements and convexity of Level sets in
PDES, Lecture notes in Mathematics, 1150, Springer-Verlag,
Heidelberg 1985.

\bibitem{ref SK}S. Kesavan, Symmetrization and Applicantions. Series in Anylysis, Vol.3,
           World Scientific Books, April,2006.

\bibitem{ref KJ1}M.-Th. Kohler-Jobin, Isoperimetric monotonicity and
Isoperimetric inequalities of Payne-Rayner type for the first
eigenfunction of the Helmholtz problem, J. Appl. Math. Phys. (ZAMP),
32(1981), p625-646.

\bibitem{ref KJ2}M.-Th. Kohler-Jobin, Sur la premi\`{e} function
propre d'une membrane: une extension \`{a} N dimensions de
l'in\'{e}galit\'{e} isop\'{e}rim\'{e}trique de Payne-Rayner, Z.
angew. Math. Phys., 28(1977), p1137-1140.

\bibitem{ref Kra} E. Krahn, \"{U}ber eine von Rayleigh formulierte
Minimaleigenschaft des Kreises, Math. Ann. 94(1925), p97-100.

\bibitem{ref Kra1} E. Krahn, \"{U}ber Minimaleigenschaften der Kugel in drei und mehr
Dimensionen, Acta Comm. Univ. Tartu (Dorpat) A9(1926), p1-44.


\bibitem{ref BM}B. Messano, Symmetrization results for classes of nonlinear
                equations with $q-$ growth in the gradient, Nonlinear Anlysis, 64(2006), p2688-2703.


\bibitem{ref PR1} L. E. Payne, M. E. Rayner, An isoperimetric inequality
                  for the first eigenfunction in the fixed membrane problem, Journal
                  of Applied Mathematics and Physics, 23(1972), p13-15.

\bibitem{ref PR2}L. E. Payne, M. E. Rayner, Some
                 isoperimetric norm bounds for solutions of the Helmholtz equation,
                 Journal of Applied Mathematics and Physics, 24(1973), p106-110.

\bibitem{ref Sze} G. Szeg\"{o}, Inequalities for certain eigenvalues of a
membrane of given area, J. Rational Mech. Anal., 3(1954), p343-356.

\bibitem{ref Tal}G. Talenti, Elliptic equations and rearrangements,
Ann. Scuola Norm. Sup. Pisa, 3(1976), p697-718.

\bibitem{ref GT}G. Talenti, Nonlinear elliptic equations, rearrangements
               of functions and orlicz spaces, Annali di Matematica Pura ed
               Applicata, 120(1979), p0373-3114.

\bibitem{ref Wein1} H. F. Weinberger, An isoperimetric inequality for the
$n$-dimensional free membrane problem, J. Rational Mech. Anal. 5
(1956), p633-636.

\bibitem{ref Wein}H. F. Weinberger, Symmetrization in uniformly elliptic problems, In "Studies in Math. Anal.,"
Stanford Uni. Press, 1962.



\end{thebibliography}
\end{document}